\documentclass[10.5pt]{article}
\usepackage{amsmath}
\usepackage{amssymb}
\usepackage{amscd}
\usepackage{xspace}
\usepackage{verbatim}
\newcommand{\mysection}[1]{
\section{#1}\setcounter{equation}{0}}
\title{\bf Boundary Harnack inequality and a priori estimates of singular solutions of quasilinear elliptic equations}
\author{ {\bf Marie-Fran\c{c}oise Bidaut-V\'eron}, {\bf Rouba Borghol}, {\bf Laurent V\'eron}\\
{\small Department of Mathematics,}\\
 {\small  University of Tours,  FRANCE}
}%%

\date{}
\begin{document}
\maketitle

%% FONT commands
\newcommand{\txt}[1]{\;\text{ #1 }\;}%% Used in math only
\newcommand{\tbf}{\textbf}%% Bold face. Usage: \tbf{...}
\newcommand{\tit}{\textit}%% Italic
\newcommand{\tsc}{\textsc}%% Small caps
\newcommand{\trm}{\textrm}
\newcommand{\mbf}{\mathbf}%% Math bold
\newcommand{\mrm}{\mathrm}%% Math Roman
\newcommand{\bsym}{\boldsymbol}%% Bold math symbol
%%Macros for changing font size in math.
\newcommand{\scs}{\scriptstyle}%% as in subscript
\newcommand{\sss}{\scriptscriptstyle}%% as in sub-subscript
\newcommand{\txts}{\textstyle}
\newcommand{\dsps}{\displaystyle}
%%Macros for changing font size in text.
\newcommand{\fnz}{\footnotesize}
\newcommand{\scz}{\scriptsize}
%%\tiny<\scz<\fsz<\small<\large<\Large<\huge<\Huge
%%%%%%%%%%%%
%%%%%%%%%%%%
%% EQUATION commands
\newcommand{\be}{
\begin{equation}
}
\newcommand{\bel}[1]{
\begin{equation}
\label{#1}}
\newcommand{\ee}{
\end{equation}
}%% This macro does not work with amstex.
\newcommand{\eqnl}[2]{
\begin{equation}
\label{#1}{#2}
\end{equation}
}%%use not advisable; confusing
%%%%%%%%%%%%%%%
%% Unnumbered THEOREM env.
%% New env. to be used for unnumbered theorem, lemma etc. (but with specified name)
\newtheorem{subn}{\name}
\renewcommand{\thesubn}{}
\newcommand{\bsn}[1]{\def\name{#1}
\begin{subn}}
\newcommand{\esn}{
\end{subn}}
%%%%%%%%%%%%%%
%% NUMBERED THEOREM env.
%% Environments: theorem, lemma, corollary defintion and related commands,
%% designed to provide consecutive numbering of these forms.
\newtheorem{sub}{\name}[section]
\newcommand{\dn}[1]{\def\name{#1}}   %used in conjuction with sub or subn.
\newcommand{\bs}{
\begin{sub}}
\newcommand{\es}{
\end{sub}}
\newcommand{\bsl}[1]{
\begin{sub}\label{#1}}
%% the above must be preceeded by \dn (name definition),
%% however this is superceded by the list of commands bth etc.  below.
%%%%%%%%%%%%
%% NUMBERED THEOREM env. (cont.)
%% List of commands derived from 'sub' env. for theorem, lemma etc.
%% designed to provide consecutive numbering of these forms.
\newcommand{\bth}[1]{\def\name{Theorem}
\begin{sub}\label{t:#1}}
\newcommand{\blemma}[1]{\def\name{Lemma}
\begin{sub}\label{l:#1}}
\newcommand{\bcor}[1]{\def\name{Corollary}
\begin{sub}\label{c:#1}}
\newcommand{\bdef}[1]{\def\name{Definition}
\begin{sub}\label{d:#1}}
\newcommand{\bprop}[1]{\def\name{Proposition}
\begin{sub}\label{p:#1}}
%%%%%%%%%%%%%%%%%%%%%%%%%%%%%%%%%%
%% RERERENCE commands.
%% \newcommand{\R}[1]{(\ref{#1})}
\newcommand{\R}{\eqref}
\newcommand{\rth}[1]{Theorem~\ref{t:#1}}
\newcommand{\rlemma}[1]{Lemma~\ref{l:#1}}
\newcommand{\rcor}[1]{Corollary~\ref{c:#1}}
\newcommand{\rdef}[1]{Definition~\ref{d:#1}}
\newcommand{\rprop}[1]{Proposition~\ref{p:#1}}
%%%%%%%%%%%
%% ARRAY commands.
\newcommand{\BA}{
\begin{array}}
\newcommand{\EA}{
\end{array}}
\newcommand{\BAN}{\renewcommand{\arraystretch}{1.2}
\setlength{\arraycolsep}{2pt}
\begin{array}}
\newcommand{\BAV}[2]{\renewcommand{\arraystretch}{#1}
\setlength{\arraycolsep}{#2}
\begin{array}}
%Note: The first variable gives the amount of stretching: (#1) x default.
%For instance #1=1.2 means a 20% stretching. The second variable should be
%written for instance in the form  4pt ; here the default is 5pt
%\newcommand{\EAN}{\end{array}\setlength{\arraycolsep}{5pt}}
\newcommand{\BSA}{
\begin{subarray}}
\newcommand{\ESA}{
\end{subarray}}
%Note: These are used in subscripts as well as superscripts. They work essentially
%% like 'array'.
\newcommand{\BAL}{
\begin{aligned}}
\newcommand{\EAL}{
\end{aligned}}
\newcommand{\BALG}{
\begin{alignat}}
\newcommand{\EALG}{
\end{alignat}}%% the abbrev. does not work with latex2e
\newcommand{\BALGN}{
\begin{alignat*}}
\newcommand{\EALGN}{
\end{alignat*}}%% the abbrev. does not work with latex2e
%% The 'aligned' environment must be placed inside an 'equation' env.
%% in the same way as the array.
%% One could use also the 'align' env. or the 'alignat' env.
%% However in this case each line is numbered, unless '\notag' is used.
%% The 'alignat'
%% has a slightly different format (the number of columns must be specified in advance)
%% but it has the advantage that the distance between columns is at our disposition.
%% (The default would be zero distance.) Using 'alignat*' we can have the advantages
%% of alignat plus the situation where separate lines are not numbered.
%% However in this case there is no numbering at all (unless we provide a tag).
%%%%%%%%%%
%% PROOF, REMARK etc.
\newcommand{\note}[1]{\textit{#1.}\hspace{2mm}}
\newcommand{\Proof}{\note{Proof}}
\newcommand{\qeda}{\hspace{10mm}\hfill $\square$}
\newcommand{\qed}{\\
${}$ \hfill $\square$}
\newcommand{\Remark}{\note{Remark}}
%%%%%%%% Style command.
\newcommand{\modin}{$\,$\\
[-4mm] \indent}
%% To be used after \mysection in order to start new line with \indent.
%%%%%%%%%%%%
%% MATHEMATICAL symbols
\newcommand{\forevery}{\quad \forall}
\newcommand{\set}[1]{\{#1\}}
\newcommand{\setdef}[2]{\{\,#1:\,#2\,\}}
\newcommand{\setm}[2]{\{\,#1\mid #2\,\}}
%% Arrows
\newcommand{\lra}{\longrightarrow}
\newcommand{\lla}{\longleftarrow}
\newcommand{\llra}{\longleftrightarrow}
\newcommand{\Lra}{\Longrightarrow}
\newcommand{\Lla}{\Longleftarrow}
\newcommand{\Llra}{\Longleftrightarrow}
\newcommand{\warrow}{\rightharpoonup}
%% Brackets, delimiters
\newcommand{
\paran}[1]{\left (#1 \right )}%% adjustable parantheses
\newcommand{\sqbr}[1]{\left [#1 \right ]}%% adjustable square brackets
\newcommand{\curlybr}[1]{\left \{#1 \right \}}%% adjustable curly brackets
\newcommand{\abs}[1]{\left |#1\right |}%% adjustable vertical delimiters
\newcommand{\norm}[1]{\left \|#1\right \|}%% adjustable norm
\newcommand{
\paranb}[1]{\big (#1 \big )}%% non-adjustable parantheses (big)
\newcommand{\lsqbrb}[1]{\big [#1 \big ]}%% non-adjustable square brackets (big)
\newcommand{\lcurlybrb}[1]{\big \{#1 \big \}}%% non-adjustable curly brackets (big)
\newcommand{\absb}[1]{\big |#1\big |}%% non-adjustable vertical delimiters (big)
\newcommand{\normb}[1]{\big \|#1\big \|}%% non-adjustable norm (big)
\newcommand{
\paranB}[1]{\Big (#1 \Big )}%% non-adjustable parantheses (Big)
\newcommand{\absB}[1]{\Big |#1\Big |}%% non-adjustable vertical delimiters (Big)
\newcommand{\normB}[1]{\Big \|#1\Big \|}%% non-adjustable norm (Big)

%%%%%%%%%%%%%%%%%
%% Adjustable parantheses etc. in a different DEFINITION format.
%\def\adp(#1){\left (#1 \right )}%% adjustable parantheses
%\def\adsb(#1){\left [#1\right ]}%% adjustable square brackets
%\def\adcb(#1){\left \{#1\right \}}%% adjustable curly brackets
%\def\abs|#1|{\left |#1\right |}%% adjustable vertical delimiters
%%%%%%%%%%%%%%%%
%% More mathematical symbols
\newcommand{\thkl}{\rule[-.5mm]{.3mm}{3mm}}
\newcommand{\thknorm}[1]{\thkl #1 \thkl\,}
\newcommand{\trinorm}[1]{|\!|\!| #1 |\!|\!|\,}
\newcommand{\bang}[1]{\langle #1 \rangle}%% angle bracket
\def\angb<#1>{\langle #1 \rangle}%% angle bracket
%% The two last lines yield the same result.
%% The second is used as follows: \angb<a,b>
\newcommand{\vstrut}[1]{\rule{0mm}{#1}}
\newcommand{\rec}[1]{\frac{1}{#1}}
%% OPERATOR names.
%% OPERATOR names.
\newcommand{\opname}[1]{\mbox{\rm #1}\,}
\newcommand{\supp}{\opname{supp}}
\newcommand{\dist}{\opname{dist}}
\newcommand{\myfrac}[2]{{\displaystyle \frac{#1}{#2} }}
\newcommand{\myint}[2]{{\displaystyle \int_{#1}^{#2}}}
\newcommand{\mysum}[2]{{\displaystyle \sum_{#1}^{#2}}}
\newcommand {\dint}{{\displaystyle \int\!\!\int}}%%%%%%%%%%
%%%%%%% SPACE commands
\newcommand{\q}{\quad}
\newcommand{\qq}{\qquad}
\newcommand{\hsp}[1]{\hspace{#1mm}}
\newcommand{\vsp}[1]{\vspace{#1mm}}
%%%%%%%%%%%
%% ABREVIATIONS
\newcommand{\ity}{\infty}
\newcommand{\prt}{
\partial}
\newcommand{\sms}{\setminus}
\newcommand{\ems}{\emptyset}
\newcommand{\ti}{\times}
\newcommand{\pr}{^\prime}
\newcommand{\ppr}{^{\prime\prime}}
\newcommand{\tl}{\tilde}
\newcommand{\sbs}{\subset}
\newcommand{\sbeq}{\subseteq}
\newcommand{\nind}{\noindent}
\newcommand{\ind}{\indent}
\newcommand{\ovl}{\overline}
\newcommand{\unl}{\underline}
\newcommand{\nin}{\not\in}
\newcommand{\pfrac}[2]{\genfrac{(}{)}{}{}{#1}{#2}}% frac with parantheses.
%%%%%%%%%%%
%%%%%%%%%%%%%

%%Macros for Greek letters.
\def\ga{\alpha}     \def\gb{\beta}       \def\gg{\gamma}
\def\gc{\chi}       \def\gd{\delta}      \def\ge{\epsilon}
\def\gth{\theta}                         \def\vge{\varepsilon}
\def\gf{\phi}       \def\vgf{\varphi}    \def\gh{\eta}
\def\gi{\iota}      \def\gk{\kappa}      \def\gl{\lambda}
\def\gm{\mu}        \def\gn{\nu}         \def\gp{\pi}
\def\vgp{\varpi}    \def\gr{\rho}        \def\vgr{\varrho}
\def\gs{\sigma}     \def\vgs{\varsigma}  \def\gt{\tau}
\def\gu{\upsilon}   \def\gv{\vartheta}   \def\gw{\omega}
\def\gx{\xi}        \def\gy{\psi}        \def\gz{\zeta}
\def\Gg{\Gamma}     \def\Gd{\Delta}      \def\Gf{\Phi}
\def\Gth{\Theta}
\def\Gl{\Lambda}    \def\Gs{\Sigma}      \def\Gp{\Pi}
\def\Gw{\Omega}     \def\Gx{\Xi}         \def\Gy{\Psi}

%%Macros for calligraphic letters.
\def\CS{{\mathcal S}}   \def\CM{{\mathcal M}}   \def\CN{{\mathcal N}}
\def\CR{{\mathcal R}}   \def\CO{{\mathcal O}}   \def\CP{{\mathcal P}}
\def\CA{{\mathcal A}}   \def\CB{{\mathcal B}}   \def\CC{{\mathcal C}}
\def\CD{{\mathcal D}}   \def\CE{{\mathcal E}}   \def\CF{{\mathcal F}}
\def\CG{{\mathcal G}}   \def\CH{{\mathcal H}}   \def\CI{{\mathcal I}}
\def\CJ{{\mathcal J}}   \def\CK{{\mathcal K}}   \def\CL{{\mathcal L}}
\def\CT{{\mathcal T}}   \def\CU{{\mathcal U}}   \def\CV{{\mathcal V}}
\def\CZ{{\mathcal Z}}   \def\CX{{\mathcal X}}   \def\CY{{\mathcal Y}}
\def\CW{{\mathcal W}} \def\CQ{{\mathcal Q}} 
%%%%%
%%Macros for 'blackboard' letters (See (27) for display.)
\def\BBA {\mathbb A}   \def\BBb {\mathbb B}    \def\BBC {\mathbb C}
\def\BBD {\mathbb D}   \def\BBE {\mathbb E}    \def\BBF {\mathbb F}
\def\BBG {\mathbb G}   \def\BBH {\mathbb H}    \def\BBI {\mathbb I}
\def\BBJ {\mathbb J}   \def\BBK {\mathbb K}    \def\BBL {\mathbb L}
\def\BBM {\mathbb M}   \def\BBN {\mathbb N}    \def\BBO {\mathbb O}
\def\BBP {\mathbb P}   \def\BBR {\mathbb R}    \def\BBS {\mathbb S}
\def\BBT {\mathbb T}   \def\BBU {\mathbb U}    \def\BBV {\mathbb V}
\def\BBW {\mathbb W}   \def\BBX {\mathbb X}    \def\BBY {\mathbb Y}
\def\BBZ {\mathbb Z}

%%Macros for Ghotic (Fraktur) letters.
\def\GTA {\mathfrak A}   \def\GTB {\mathfrak B}    \def\GTC {\mathfrak C}
\def\GTD {\mathfrak D}   \def\GTE {\mathfrak E}    \def\GTF {\mathfrak F}
\def\GTG {\mathfrak G}   \def\GTH {\mathfrak H}    \def\GTI {\mathfrak I}
\def\GTJ {\mathfrak J}   \def\GTK {\mathfrak K}    \def\GTL {\mathfrak L}
\def\GTM {\mathfrak M}   \def\GTN {\mathfrak N}    \def\GTO {\mathfrak O}
\def\GTP {\mathfrak P}   \def\GTR {\mathfrak R}    \def\GTS {\mathfrak S}
\def\GTT {\mathfrak T}   \def\GTU {\mathfrak U}    \def\GTV {\mathfrak V}
\def\GTW {\mathfrak W}   \def\GTX {\mathfrak X}    \def\GTY {\mathfrak Y}
\def\GTZ {\mathfrak Z}   \def\GTQ {\mathfrak Q}

\font\Sym= msam10 % special symbols
\def\SYM#1{\hbox{\Sym #1}}
\newcommand{\bdw}{\prt\Gw\xspace}
\medskip
\mysection {Introduction}
Let $\Gw$ be a domain is $\BBR^N$, $N\geq 2$, $d$ a locally bounded and measurable function defined in $\Gw$ and $p$ a real number larger than $1$. This article deals with the study of positive solutions of
\begin {equation}\label {p-lapl}
-div\left(\abs {Du}^{p-2}Du\right)-d(x)u^{p-1}=0\quad \mbox {in }\Gw
\end {equation}
which admit an isolated  singularity on the boundary of $\Gw$. It is known since the starting pioneering work of Serrin \cite {Se1} that one of the
main goals for studying the regularity of solutions of quasilinear equations consists in obtaining Harnack inequalities. The simplest form of this inequality is the following: Assume $B_{2r}\subset\subset \Gw$ and
$d\in L^{\infty}(B_{2r})$, then there exists a constant 
$C=C(N,p,r\left({\norm {d}_{ L^{\infty}(B_{2r})}}\right)^{1/p})\geq 1$ such that any nonnegative solution $u$ of (\ref {p-lapl}) in $B_{2r}$ satisfies 
\begin {equation}\label {harn1}
u(x)\leq Cu(y)\forevery (x,y)\in B_r\ti  B_r.
\end {equation}
Actually this inequality is valid for a much wider class of operators in divergence form with a power-type growth. Among the important consequences of this inequality are the H\"older continuity of the weak solutions of (\ref {p-lapl}) and the two-side estimate of solutions admitting an isolated singularity. Among more sophisticated consequences are the obtention of local upper estimates of solutions of the same equation near a singular point. This program has been carried out by Gidas and Spruck 
for equation
\begin {equation}\label {semilin}
-\Gd u=u^q
\end {equation}
in the case $N\geq 2$ and $1<q<(N+2)/(N-2)$ \cite {GS}, and recently by Serrin and Zhou \cite {Se2} for equation 
\begin {equation}\label {p-lapl2}
-div\left(\abs {Du}^{p-2}Du\right)=u^q
\end {equation}
in the case $N\geq p$ and  $p-1<q<Np/(N-p)-1$. A third type of applications of Harnack inequality linked to the notion of isotropy leads to the description of positive isolated singularities of solutions. This was carried out by V\'eron \cite {Ve1} for 
\begin {equation}\label {semilin2}
-\Gd u+u^q=0
\end {equation}
in the case $1<q<N/(N-2)$, and by Friedman and V\'eron \cite {FV} 
for
\begin {equation}\label {p-lapl3}
-div\left(\abs {Du}^{p-2}Du\right)+u^q=0
\end {equation}
when  $p-1<q<N(p-1)/(N-p)$. When the singularity of $u$ is not an internal point but a boundary point, the situation is more complicated and the mere inequality (\ref {harn1}) with only one function has no meaning. Boundary Harnack inequalities which deals with two nonnegative solutions of (\ref {p-lapl}) vanishing on a part of the boundary asserts that the two solutions must vanish at the same rate. For linear second order elliptic equations they are used for studying the properties of the harmonic measure \cite {CFMS} (see also \cite {Ba}). For p-harmonic function in a ball, a sketch of construction is given by Manfredi and Weitsman \cite {MW} in order to obtain Fatou type results. In this article we consider singular solutions of (\ref {p-lapl}) with a singular potential type reaction term. The first result we prove is the following: {\it Assume $\prt\Gw$ is $C^2$ and $d$ is measurable,  locally bounded in  $\overline \Gw\setminus\{a\}$ for some $a\in\prt\Gw$ and satisfies
\begin {equation}\label {d-est}
\abs {d(x)}\leq C_0\abs {x-a} ^{-p}\quad \mbox {a. e. in }\;B_R(a)\cap \Gw
\end {equation}
for some $C_{0}, R>0$ and $p>1$. Then there exists a positive constant
$C$ depending also on $N$, $p$ and $C_0$ such that if 
$u\in C^1(\overline\Gw\setminus\{a\})$ is a  
nonnegative solution of (\ref {p-lapl}) vanishing on $\prt\Gw\setminus \{a\}$, there holds
\begin {equation}\label {bdrH1}
\myfrac{u(y)}{C\gr(y)}\leq \myfrac {u(x)}{\gr(x)}\leq 
\myfrac{Cu(y)}{\gr(y)}
\forevery (x,y)\in \Gw\ti\Gw\quad {s. t. }\;\abs x=\abs y.
\end {equation}
where $\gr(.)$ is the distance function to $\prt\Gw$}. Another form of this estimate, usually called boundary Harnack inequality, asserts that if $u_{1}$ and $u_{2}$ are two nonnegative solutions of (\ref {p-lapl}) vanishing on $\prt\Gw\setminus \{a\}$, there holds
\begin {equation}\label {bdrH2}
\myfrac{u_{1}(x)}{Cu_{1}(y)}\leq \myfrac{u_{2}(x)}{u_{2}(y)}\leq 
\myfrac{Cu_{1}(x)}{u_{1}(y)}
\forevery (x,y)\in \Gw\ti\Gw\quad {s. t. }\;\abs x=\abs y.
\end {equation}
for some structural constant $C>0$. Another consequence of the construction leading to (\ref {bdrH1}) is the existence of a power-like a priori estimate: {\it Assume $\Gw$ is a bounded $C^2$ domain with $a\in \prt\Gw$, $A\in\Gw$ is an arbitrary point and $d$ a measurable function such that
 \begin {equation}\label {d-est1}
\abs {d(x)}\leq C_0\abs {x-a} ^{-p}\quad \mbox {a. e. in }\; 
\Gw. 
\end {equation}
Then there exist two positive constants
$\ga>0$ depending on $N$, $p$, $\Gw$ and $C_0$, and $C$ depending on the same parameters and also on $A$ such that, any nonnegative solution $u\in C(\overline\Gw\setminus\{a\})\cap W^{1,p}_{loc}(\overline\Gw\setminus\{a\}))$ which vanishes on
$\prt\Gw\setminus \{a\}$  verifies}
\begin {equation}\label {bdrH3}
u(x)\leq C\myfrac {\gr(x)}{\abs {x-a}^{\ga+1}}u(A)\forevery x\in
\overline\Gw\setminus\{a\}.
\end {equation}

The precise value of $\ga$ is unknown and surely difficult to know explicitely, even in the simplest case when $u$ is a p-harmonic function. In several cases the value of $\ga$ is associated to the construction of separable $p$-harmonic functions called the spherical $p$-harmonics. Another striking applications of the boundary Harnack principle deals with the structure of the set of positive singular solutions. We prove the following: {\it
Let $\Gw$ be $C^2$ and bounded, $a\in\prt\Gw$ and $d$ satisfies (\ref {d-est1}). Assume also that the operator 
$v\mapsto -div\left(\abs {Dv}^{p-2}Dv\right)-d(x)v^{p-1}$ satisfies the comparison principle in
$\Gw\setminus B_{\ge}(a)$ for any $\ge>0$, among nonnegative solutions which vanishes on $\prt\Gw\setminus B_{\ge}(a)$. If $u$ and $v$ are two positive solutions of (\ref {p-lapl}) in $\Gw$ which vanish on
$\prt\Gw\setminus\{a\}$, there exists $k>0$ such that 
\begin {equation}\label {bdrH4}
k^{-1}u(x)\leq v(x)\leq ku(x)\forevery x\in\Gw.
\end {equation}
Furthermore, if we assume also either $p=2$, either $p>2$ and $u$ has no critical point in $\Gw$, or 
$1<p<2$ and $d\geq 0$,  there exists $k>0$ such that }
\begin {equation}\label {bdrH5}
v(x)= ku(x)\forevery x\in\Gw.
\end {equation}
In the last section we give some partial results concerning the existence of singular solutions of equations of type 
(\ref {p-lapl}) and their link with separable solutions which are solution under the form $u(x)=\abs x^{-\gg}\phi(x/\abs x)$. If $d\equiv 0$ such specific solutions, studied by Kroll and Maz'ya \cite {KM}, Tolksdorff \cite {To}, Kichenassamy and V\'eron \cite {KV}, are called spherical $p$-harmonics.

Our paper is organized as follows: 1- Introduction. 2 The boundary Harnack principle. 3 Existence of singular solutions. 4 References.

\mysection {The boundary Harnack principle}
In this section we consider nonnegative solutions of
\begin {equation}\label {main}
-div\left(\abs {Du}^{p-2}Du\right)-d(x)u^{p-1}=0
\end {equation}
in a domain  $\Gw$ which may be Lipschitz continuous or $C^2$. The function $d$, is supposed to be measurable and singular in the sense that it satisfies 
\begin {equation}\label {main-d}
\abs {d(x)}\leq C_0\abs {x-a}^{-p}\quad\mbox {a.e. }\; \in \Gw
\end {equation}
for some point $a\in\prt \Gw$ and some $C_0\geq 0$. By a solution of (\ref {main}) vanishing on $\prt\Gw\setminus \{a\}$, we mean a 
$u\in C(\bar\Gw\setminus\{a\})$ such that $Du\in L^p(K)$ for every $K$ compact, $K\subset \bar\Gw\setminus\{a\}$ which verifies
\begin {equation}\label {main-test}
\int_{\Gw}\left(\abs {Du}^{p-2}Du.D\gz+d(x)u^{p-1}\gz\right)\, dx=0
\end {equation}
for every $\gz\in C^1(\bar\Gw)$, with compact support in $\bar\Gw\setminus\{a\}$.

\subsection {Estimates near the boundary in Lipschitz domains}
Let  $\Gw$ be a bounded domain in $\BBR^N$ with a Lipschitz continuous boundary. Then there exist $m>1$ and $r_0>0$ such that for any $Q\in\prt\Gw$ there exists an isometry $\CI_Q$ in $\BBR^N$ and a Lipschitz continuous real valued function $\phi$ defined in $\BBR^{N-1}$ such that 
\begin {equation}\label {lip1}
\abs {\phi(x)-\phi(y)}\leq m\abs {x-y}\forevery (x,y)\in 
\BBR^{N-1}\ti \BBR^{N-1},
\end {equation}
and 
$$B_{2r_0}(Q)\cap\{x=(x',x_N)=(x_1,...,x_N):x_N\geq \phi(x')
\}=\CI_Q\left(\Gw\cap B_{2r_0}(Q)\right).
$$
For any $A\in B_{r_{0}(Q)}\cap\prt\Gw$, $A=(a',\phi(a')$, $r>0$ and $\gg>0$, we denote by $C_{A,r,\gg}$
 the truncated cone
$$
C_{A,r,\gg}=\{x=(x',x_{N}): x_{N}>\phi(a')\,,\; \gg\abs {x'-a'}\leq x_{N}-\phi(a')\}\cap B_{r}(A)
$$
The opening angle of this cone is $\gth_{\gg}=\tan^{-1}(1/\gg)$
Clearly, for every  $\gg\geq m$ and $0<r\leq r_{0}$, $\CI^{-1}_Q(C_{A,r,\gg})$ is included into $\Gw$. Up to an orthogonal change of variable, we shall assume that $\CI_Q=Id$. We denote also by $\rho(x)$ the distance from $x$ to $\prt\Gw$. The next result is a standard geometric construction which can be found in \cite {Ba}
% %%LEMMA   "geom1"  %%CONNECTED CHAINS%%%%%%%%%%%%%%%%%%%%
%%%%%%%%%%%%%%%%%%%%%%%%%%%%%%%%%%%%%%%%%%%
\blemma{geom1} Let $Q\in\prt\Gw$ and $0<r<r_0/5$ and $h>1$ an integer. There exists $N_0\in\BBN$ depending only on $m$ such that for any points $x$ and $y$ in $\Gw\cap B_{3r/2}(Q)$ verifying $\min\{\gr(x),\gr(y)\}\geq r/2^{h}$, there exists a connected chain of balls $B_1,...B_j$ with $j\leq N_0h$ such that 
\begin {equation}\label {incl}
x\in B_1,\; y\in B_j,\;B_i\cap B_{i+1}\neq \emptyset  \mbox { and } 2B_i\subset B_{2r}(Q)\cap\Gw \mbox { for }1\leq i\leq j-1.
\end {equation}
\es
%%%PROOF  "geom1"  %%CONNECTED CHAINS%%%%%%%%%%%%%%%%%%%%
%%%%%%%%%%%%%%%%%%%%%%%%%%%%%%%%%%%%%%%%%%%

%%%LEMMA   "holdest"  %%HARNACK INTERIOR%%%%%%%%%%%%%%%%%%%%
%%%%%%%%%%%%%%%%%%%%%%%%%%%%%%%%%%%%%%%%%%%
\blemma {holdest} Assume $u$ is a nonnegative solution of (\ref{main}) in $B_{2r}(P)$ where
$\abs {a-P}\geq 4r$. Then there exists a positive constant $c_1$ depending on $p$, $N$ and $C_0$
such that 
\begin {equation}\label {harn-int}
u(x)\leq c_1u(y)\forevery (x,y)\in B_{ r}(P)\ti B_{ r}(P).
\end {equation}
\es
%%%PROOF%%HARNACK INTERIOR%%%%%%%%%%%%%%%%%%%%
%%%%%%%%%%%%%%%%%%%%%%%%%%%%%%%%%%%%%%%%%%%
\Proof By a result of Trudinger \cite [Th. 1.1]{Tr}, if $u$ is a nonnegative solution of (\ref {main}) in 
$B_{2r}(P)$ there exists a constant $C'$ depending on $N$, $p$, $r$ and $\norm {\abs d^{1/p}}_{L^{\ity}(B_{2r}(P))}$ such that (\ref {harn-int}) holds.
Furthermore $C'\leq C_2\exp \left({C_3r\norm {\abs d^{1/p}}_{L^{\ity}(B_{2r}(P))}}\right)$ where
$C_2$ and $C_3$ depend on $N$ and $p$. This implies (\ref {harn-int}) since, by (\ref {main-d}), $r\norm {\abs d^{1/p}}_{L^{\ity}(B_{2r}(P))}$ remains bounded by a constant depending on $p$ and $C_0$. \qeda \\

%%%LEMMA  "holdestII"%%HARNACK INTERIOR II%%%%%%%%%%%%%%%%%%%%
%%%%%%%%%%%%%%%%%%%%%%%%%%%%%%%%%%%%%%%%%%%
Up to a translation, we shall assume that the singular boundary point $a$ is the origin of coordinates.\\

\blemma {holdestII} Assume $\Gw$ is as in \rlemma {geom1} with $Q\neq 0\in\prt\Gw$,
$0<r\leq \min\{r_0,\abs{Q}/4\}$ and  $u$ is a nonnegative solution of (\ref{main}) in $B_{2r}(Q)$. Then there exists a positive constant $c_2>1$ depending on $p$, $N$, $m$ and $C_0$
such that 
\begin {equation}\label {harn-h}
u(x)\leq c_2^hu(y),
\end {equation}
for every $x$ and $y$ in $B_{3r/2}(Q)\cap\Gw$ such that $\min\{\gr(x),\gr(y)\}\geq r/2^h$ for some 
$h\in\BBN$.
\es
\Proof By \rlemma {geom1} there exists $N_0\in\BBN^{*}$ and a connected chain of $j\leq N_0h$ balls $B_i$ with respective radii $r_i$ and centers $x_i$, satisfying (\ref {incl}). Thus
\begin {equation}\label {harn-i}
\max_{B_i} u\leq c_1\min_{B_i} u\forevery i=1,...,j,
\end {equation}
by the previous lemma. Therefore (\ref {harn-h}) holds with $c_2=c_1^{N_0}$.\qeda \\

%%%LEMMA  "holdestIII"%%HARNACK-HOLDER%%%%%%%%%%%%%%%%%%%%
%%%%%%%%%%%%%%%%%%%%%%%%%%%%%%%%%%%%%%%%%%%
\blemma {holdestIII} Let $0<r\leq \abs{Q}/4$ and $u$ be a nonnegative solution of 
(\ref {main}) in $B_{2r}(Q)\cap\Gw$ which vanishes on $\prt\Gw\cap B_{2r}(Q)$. If $P\in \prt\Gw\cap B_{r}(Q)$ and $0<s<r/(1+m)$ so that $B_s(P)\subset B_{2r}(Q)$, there exist two positive constants 
$\gd$ and $c_3$ depending on $N$, $p$, $m$ and $C_0$ such that
\begin {equation}\label {harn-hold}
u(x)\leq c_3\frac {\abs {x-P}^\gd}{s^\gd}M_{s,P}(u)
\end {equation}
for all $x\in B_s(P)\cap\Gw$, where $M_{s,P}(u)=\max\{u(z):z\in B_s(P)\cap\Gw\}$. \es
%%%%%%PROOF  "holdestIII"%%HARNACK-HOLDER%%%%%%%%%%%%%%%%%%%%
% %%%%%%%%%%%%%%%%%%%%%%%%%%%%%%%%%%%%%%%%%
\Proof Since $\prt \Gw$ is Lipschitz, it is regular in the sense that there exists $\gth>0$, $s_{1}>0$ such that
$$ meas\,(\Gw^c\cap B_{s}(y))\geq\gth\, meas\,( B_{s}(y)), \forevery y\in\prt\Gw,\;\forall 0<s<s_{1}.
$$
By \cite[Th. 4.2]{Tr} there exists $\gd\in (0,1)$ depending on  $p$, $N$, $C_{0}$, $\gth$ and $s_{1}$, such that for any $y\in\prt\Gw$, there holds
\begin {equation}\label {osc}
\abs {u(z)-u(z')}\leq C\left(\myfrac {s}{s_{1}}\right)^{\gd(1-\gg)}\forevery (z,z')
\in B_{s}(y)\cap\Gw\ti B_{s}(y)\cap\Gw,
\end {equation}
where $C$ depends on $p$, $N$, $C_{0}$ and $\sup_{B_{s_{1}}(y)\cap \Gw} u =M_{s_{1},y}(u)$. Because the equation is homogeneous with respect to $u$, this local estimate is invariant if we replace $u$ by
$\tilde u=u/M_{s_{1},y}(u) $. Thus the dependence is homogeneous of degree 1 with respect to $M_{s_{1},y}(u) $, which implies
\begin {equation}\label {osc-2}
\abs {u(z)-u(z')}\leq C'\left(\myfrac {s}{s_{1}}\right)^{\gd(1-\gg)}M_{s_{1},y}(u)\forevery (z,z')
\in B_{s}(y)\cap\Gw\ti B_{s}(y)\cap\Gw.
\end {equation}
Taking $z'=P=y$, $s=\abs {x-P}$, we derive (\ref {harn-hold}).\qeda\\

If $X\in B_{2r_{0}}(0)\cap\prt\Gw$ and $r>0$, we denote by $A_{r}(X)$ the point with coordinates
$(x',\phi(x')+r)$. The next result is the key point in the construction. Although it follows \cite {Ba}, we give the proof for the sake of completeness.
%%%%%%%%%%%%%%%%%%%%%%%%%%%%%%%%%%%%%%%%%%%%
%%%LEMMA  "norm-est"%%NORMAL_{ESTIMATE}%%%%%%%%%%%%%%%%%%%%

\blemma{norm} Let $0<r\leq \min\{2r_{0},\abs{Q}/8,s_{1},2^5\}$ and $u$ be a nonnegative solution of 
(\ref {main}) in $B_{2r}(Q)\cap\Gw$ which vanishes on $\prt\Gw\cap B_{2r}(Q)$. Then there exists a positive constant $c_{4}$ depending only on $N$, $p$, $m$ and $C_{0}$ such that  
\begin {equation}\label {norm-est}
u (x)\leq c_4 u(A_{r/2}(Q))\forevery x\in B_{r}(Q)\cap\Gw.
\end {equation}
\es
%%%%%%PROOF  "NORM"%%HARNACK-HOLDER%%%%%%%%%%%%%%%%%%%%
% %%%%%%%%%%%%%%%%%%%%%%%%%%%%%%%%%%%%%%%%%
\Proof  The proof is by contradiction. We first notice from (\ref{harn-hold}) that if $P\in B_{2r}(Q)\cap\prt\Gw$ verifies $B_{s}(P)\cap\Gw\subset B_{2r}(Q)\cap\Gw$ and if
$c_{5}=(2c_{3})^{1/\gd}$, there holds
\begin {equation}\label {norm-est1}
M_{s/c_{5},P}(u)\leq \myfrac {1}{2}M_{s,P}(u).
\end {equation}
By \rlemma {holdestII}, if $y\in B_{3r/2}(Q)$ satisfies $u(y)>c^h_{2}u_{A_{r/2}(Q)}$, it means that
$\gr(y)<r/2^h$. Let $M>0$ such that $2^M>c_{2}$ (defined in \rlemma {holdestII}), $N=\max \{1+\BBE(6+M\ln c_{5}/\ln 2), M+5\}$, so that $2^N>2^6c_{5}^M$, and
$c_{4}=c_{2}^{N}$. Let $u$ be a positive solution of (\ref {main}) vanishing on $B_{2r}(Q)\cap\prt\Gw$ which satisfies
\begin {equation}\label {norm-est2'}
u(Y_{0})> c_{2}^{N} u(A_{r/2}(Q)),
\end {equation}
for some $Y_{0}\in  B_{r}(Q)\cap\Gw$. Then $\gr(Y_{0})<r/2^N$. Let $Q_{0}\in\prt\Gw$ such that $\gr(Y_{0})=\abs {Y_{0}-Q_{0}}$. Then
$$\abs {Q-Q_{0}}\leq \abs {Y_{0}-Q_{0}}+\abs {Y_{0}-Q}\leq r/2^N+r\leq r(1+2^{-5})
$$
Therefore $Q_{0}\in B_{3r/2}(Q)\cap\prt\Gw$. Set $s_{2}=c^M_{5}r/2^{N}$, then $B_{s_{2}}(Q_{0})\subset B_{r(1+2^{-5}+2^{-6})}(Q)\subset
B_{3r/2}(Q)$ because $s_{2}\leq 2^{-6}$ by the choice of $N$. Applying (\ref {norm-est1}) with $s=s_{2}$
 yields to
$$M_{s_{2},Q_{0}}(u)\geq 2^M M_{s_{2}/c^M_{5},Q_{0}}(u)\geq 2^Mu(Y_{0})\geq 
2^Mc_{2}^{N} u(A_{r/2}(Q))\geq c_{2}^{N+1} u(A_{r/2}(Q)),
$$
since $\abs {Y_{0}-Q_{0}}\leq r/2^N=s_{2}/c_{5}^M$ and $2^M>c_{2}$. Hence we can choose 
$Y_{1}\in \overline {B_{s_{2}}(Q_{0})\cap\Gw}$ which realizes $M_{s_{2},Q_{0}}(u)$ and this implies that $\gr(Y_{1})<r/2^{N+1}$. A point $Q_{1}\in\prt\Gw$ such that $\gr(Y_{1})=\abs {Y_{1}-Q_{1}}$ satisfies also
$$\abs {Q-Q_{1}}\leq \abs {Q-Q_{0}}+\abs {Q_{0}-Q_{1}}\leq  r(1+2^{-5}+2^{-6}).
$$
Now
$$M_{s_{2}/2,Q_{1}}(u)\geq 2^M M_{r/2^{N+1},Q_{1}}(u)\geq 2^Mu(Y_{1})
\geq 2^{2M}c_{2}^{N} u(A_{r/2}(Q))\geq 2^{N+2}A_{r/2}(Q).
$$
Iterating this procedure, we construct two sequences $\{Y_{k}\}$ of points such that 
$\gr(Y_{k})<r/2^{k+N}$ and $\{Q_{k}\}$ such that $Q_{k}\in\prt\Gw$ and 
$\abs {Q-Q_{k}}\leq  r(1+2(2^{-5}+2^{-6}+...+2^{-5-k}))< 3r/2$ and 
$$u(Y_{k})\geq 2^{N+k}A_{r/2}(Q)\forevery k\in\BBN^{*}.
$$
Since $Y_{k}\in B_{3r/2}$ and $\gr(Y_{k})\to 0$ as $k\to\infty$ we get a contradiction with the fact that 
$A_{r/2}(Q)>0$.\qeda\\

%%%%%%%%%%%%%%%%%%%%%%%%%%%%%%%%%%%%%%%%%%%%
%%%REMARK%%%%%%%%%%%%%%}%%%%%%%%%%%%%%%%%%%%
\noindent\Remark The proof of the previous lemma shows that estimate (\ref {norm-est}) is valid for a much more general class of equations under the form
\begin {equation}\label {div}
-div A(x,u,Du)+B(x,u,Du)=0
\end {equation}
where $A$ and $B$ are respectively vector and real valued Caratheodory functions defined on $\Gw\ti\BBR\ti\BBR^N$ and verifying, for some constants $\gg>0$ and $a_{0},\,a_{1},\,C_{0}\geq 0$,
$$A(x,r,q).q\geq \gg\abs q^p,
$$
$$\abs {A(x,r,q)}\leq a_{0}\abs q^{p-1}+a_{1}\abs r^{p-1},
$$
and
$$\abs {B(x,r,q)}\leq C_{0}\abs r^{p-1}\abs {x}^{-p},
$$
for $(x,r,q)\in \Gw\ti\BBR\ti\BBR^N$.
%%%%%%%%%%%%%%%%%%%%%%%%%%%%%%%%%%%%%%%%%%%%
%%%NEW SUBSECTION %%%%%%%%%%%%%%%%%%%%%%%%%%%%%%
\subsection {Estimates near the boundary in $C^2$ domains}
From now we assume that $\Gw$ is a bounded domain with a $C^2$ boundary. For any $x\in \prt\Gw$, we denote by $\gn_{x}$ the normal unit outward vector to $\prt\Gw$
at $x$. Let $R_{0}>0$ be such that for any $x\in \prt\Gw$, the two balls 
$B_{R_{0}}(x-R_{0} \gn_{x})$ and $B_{R_{0}}(x+R_{0}\gn_{x})$ are subsets of $\Gw$ and $\bar\Gw^c$ respectively. If $P\in\prt\Gw$, we denote by $N_{r}(P)$ and 
$\CN_{r}(P)$ the points $P-r\gn_{P}$ and $=P+r\gn_{P}$. Notice that 
$r\leq R_{0}$ implies $\gr (N_{r}(P))=\gr(\CN_{r}(P))=r$.
%%%%%%%%%%%%%%%%%%%%%%%%%%%%%%%%%%%%%%%%%%%%
%%%LEMMA  TWO SIDE %%%%%%%%%%%%%%%%%%%%%%%%%%%%%%

\blemma {two-side} Let $Q\in\prt\Gw\setminus\{0\}$, $0<r\leq \min\{R_{0}/2,\abs Q/2\}$ and $u$ be a nonnegative solution of (\ref {main})\\  in $B_{2r}(Q)\cap\Gw$ which vanishes on $B_{2r}(Q)\cap\prt\Gw$. Then there exist $b\in (0,2/3)$ and $c_{6}>0$ depending respectively on $N$, $p$ and $C_{0}$ and $N$, $p$, $R_{0}$ and $C_{0}$ such that 
\begin {equation}\label {norm-est2}
\myfrac {t}{c_{6}r}\leq \myfrac {u(N_{t}(P))}{u(N_{r/2}(Q))}\leq\myfrac {c_{6}t}{r}
\end {equation}
for any $P\in B_{r}(Q)\cap\prt\Gw$ and $0\leq t\leq rb/2$.
\es
%%%%%%%%%%%%%%%%%%%%%%%%%%%%%%%%%%%%%%%%%%%%
%%%PROOF OF LEMMA  TWO SIDE %%%%%%%%%%%%%%%%%%%%%%%%%%%%%%
\Proof Up to a dilation, we can assume that $\abs Q=1$, since if we replace $x$ by 
$ x/\abs Q$, equation (\ref {main}) and the estimates (\ref {norm-est2}) are structuraly invariant 
(which means that $C_{0}$ and the $C_{i}$ are unchanged), while the curvature constant $R_{0}$ is replaced by $R_{0}/\abs Q$ which is no harm since $\Gw$ is bounded.

 \medskip

\noindent {\it Step 1} The lower bound. For $a>0$ and $\ga>0$  to be made precise later on, let us define 
$$ v(x)=V(\abs {x-N_{r/2}(P)})=
\myfrac {e^{-a(\abs {x-N_{r/2}(P)}/r)^\ga}-e^{-a/2^\ga}}
{e^{-a/4^\ga}-e^{-a/2^\ga}}
$$
for $x\in B_{r/2}(N_{r/2}(P))\cap B_{r/4}(P)$. We set $s=\abs {x-N_{r/2}(P)}$. Since $\abs Q=1$, the function $d$ satisfies $-\tilde C_{0}\leq d(x)\leq \tilde C_{0}$. Next
$$-div(\abs{Dv}^{p-2}Dv)+\tilde C_{0}v^{p-1}
=-\abs {V'}^{p-2}\left((p-1)V''+(N-1)V'/s\right)+\tilde C_{0}V^{p-1}.
$$
Therefore this last expression will be nonpositive if and only if
\begin {equation}\label {norm-est3}
(p-1)\left(\myfrac{a\ga s^{\ga}}{r^{\ga}}+1-\ga\right)+1-N
\geq \tilde C_{0}\left(\myfrac{a\ga }{r^{\ga}}\right)^{1-p}e^{(p-1)a(s/r)^{\ga}}s^{\gth}
\left(e^{-a(s/r)^{\ga}}-e^{-a/2^{\ga}}\right)^{p-1}
\end {equation}
where $\gth=p+(1-p)\ga$. But  $\gth=0$ by choosing $\ga=p/(p-1)$, thus (\ref{norm-est3})
is equivalent to
$$
(p-1)\left(\myfrac{a\ga s^{\ga}}{r^{\ga}}+1-\ga\right)+1-N
\geq\tilde C_{0}\left(\myfrac{a\ga }{r^{\ga}}\right)^{1-p}
\left(1-e^{a(1/4^{\ga}-1/2^{\ga})}\right)^{p-1}.
$$
If $r/4\leq s\leq r/2\leq 1/4$,
$$(p-1)\left(\myfrac{a\ga s^{\ga}}{r^{\ga}}+1-\ga\right)+1-N\geq 
(p-1)\left(\myfrac{a\ga}{4^{\ga}}+1-\ga\right)+1-N=\myfrac{pa}{4^{\ga}}-N,
$$
while
$$\left(\myfrac{a\ga }{r^{\ga}}\right)^{1-p}
\left(1-e^{a(1/4^{\ga}-1/2^{\ga})}\right)^{p-1}\leq 
\left(\myfrac{a\ga }{r^{\ga}}\right)^{1-p}\leq
a^{1-p}.
$$
Therefore, if we choose $a$ such that
\begin {equation}\label {norm-est4}
a^{p-1}(\myfrac{ap}{4^{\ga}}-N)\geq \tilde C_{0}, 
\end {equation}
we derive
\begin {equation}\label {norm-est5}
-div(\abs{Dv}^{p-2}Dv)+\tilde C_{0}v^{p-1}\leq 0
\end {equation}
in $B_{r/2}(N_{r/2}(P))\cap B_{r/4}(P)$. Furthermore $\gr(x)\geq r/16$ for any 
$x\in \prt B_{r/4}(P)\cap B_{r/2}(N_{r/2}(P))$, therefore
\begin {equation}\label {norm-est6}
u(x)\geq c_{2}^{-4}u(N_{r/2}(P))v(x),
\end {equation}
by \rlemma {holdestII} and since $v\leq 1$. Because $u$ is a supersolution for 
(\ref {norm-est5}), we obtain that (\ref {norm-est6}) holds for any $x\in B_{r/4}(P)\cap B_{r/2}(N_{r/2}(P))$. Finally
$$v(x)\geq \myfrac {e^{-a/2^\ga}(2^{-\ga}-(s/r)^{\ga})}{e^{-a/4^\ga}-e^{-a/2^\ga}}
\geq C(a,\ga)(1-(1-2t/r)^\ga)\geq \myfrac {C'(a,\ga)t}{r}
$$
if $x=N_{t}(P)$ with  $0\leq t\leq r/2$. This gives the left-hand side of (\ref {norm-est2}).
 \medskip

\noindent {\it Step 2} The upper bound. Let $b\in (0,2/3]$ be a parameter to be made precise later on. By the exterior sphere condition, $B_{3br}(\CN_{3rb}(P))\subset 
\bar\Gw^c$. Let $\gf_{1}$ be the first eigenfunction of the $p$-Laplace operator in 
$B_{3}\setminus\bar B_{1}$ with Dirichlet boundary conditions and $\gl_{1}$ the corresponding eigenvalue. It is well known that $\phi_{1}$ is radial. We normalize $\phi_{1}$ by 
$\phi_{1}(y)=1$ on $\{y:\abs y=2\}$ (notice that $\phi_{1}$ is radial)  and set
$$\phi_{rb}(x)=\phi_{1}\left(\myfrac{\abs {x-\CN_{rb}(P)}}{rb}\right),
$$
thus
$$-div\left(\abs {D\phi_{rb}}^{p-2}D\phi_{rb}\right)=\myfrac{\gl_{1}}{(rb)^{p}}\phi^{p-1}_{rb}
$$
in $B_{3rb}(\CN_{rb}(P))\setminus\bar B_{rb}(\CN_{rb}(P))$ and vanishes on the boundary of this domain. For $b$ small enough $\gl_{1}/(rb)^{p}\geq 1+\tilde C_{0}$ for any 
$r\in (0,1/2]$, thus
\begin {equation}\label {norm-est7}
-div\left(\abs {D\phi_{rb}}^{p-2}D\phi_{rb}\right)-\tilde C_{0}\phi_{rb}^{p-1}\geq \phi_{rb}^{p-1}
\end {equation}
in $\Gw\cap B_{3rb}(\CN_{rb}(P))\setminus\bar B_{rb}(\CN_{rb}(P))\supseteq\Gw\cap B_{2rb}(\CN_{rb}(P))$
while $u$ verifies
\begin {equation}\label {norm-est8}
-div\left(\abs {Du}^{p-2}Du\right)-\tilde C_{0}u^{p-1}\leq 0
\end {equation}
in the same domain. We can also take $b>0$ such that $B_{2br}(\CN_{br}(P))\subset B_{r}(Q)$, thus
$$u(x)\leq c_4 u(N_{r/2}(Q))
$$
for $x\in \prt B_{2rb}(\CN_{rb}(P))\cap \Gw$ by \rlemma {norm}. Now the function 
$\tilde \phi_{rb}=c_4 u(N_{r/2}(Q))\phi_{rb}$ satisfies (\ref {norm-est7}\!) in 
$\Gw\cap B_{2rb}(\CN_{rb}(P))$ and dominates $u$ on $\prt(\Gw\cap B_{2rb}(\CN_{rb}(P)))
=(\prt B_{2rb}(\CN_{rb}(P))\cap\Gw)\cup (B_{2rb}(\CN_{rb}(P))\cap\prt\Gw)$. By the Diaz-Saa inequality \cite {DS}
$$\int_{\Gw\cap B_{2rb}(\CN_{rb}(P))}\left(
\myfrac {div\left(\abs {Du}^{p-2}Du\right) }{u^{p-1}}-\myfrac {
div\left(\abs {D\tilde\phi_{rb}}^{p-2}D\tilde\phi_{rb}\right) }{\tilde\phi_{rb}^{p-1}}\right)
(u^p-\tilde\phi_{rb}^p)_{+}dx\leq 0,
$$
valid because $(u^p-\phi_{rb}^p)_{+}\in W_{0}^{1,p}(\Gw\cap B_{2rb}(\CN_{rb}(P)))$. Therefore
$$\int_{\Gw\cap B_{2rb}(\CN_{rb}(P))}(u^p-\tilde\phi_{rb}^p)_{+}dx\leq 0,
$$
from which follows the inequality  $u\leq \tilde\phi_{rb}$ in $\Gw\cap B_{2rb}(\CN_{rb}(P))$. In particular
$$u(N_{t}(P))\leq c_{4}\phi_{1}\left(\myfrac {\abs {N_{t}(P)-\CN_{rb}(P)}}{rb}\right)u(N_{r/2}(Q)).
$$
Since $\phi_{1}(s)\leq C(s-1)$ for $s\in [1,2]$, we obtain the right-hand side of 
(\ref{norm-est2}).\qeda\\

The main result of this section is the following
%%%%%%%%%%%%%%%%%%%%%%%%%%%%%%%%%%%%%%%%%%%%%%%%%%%%%%%%%%%%%%%%%%%%%%
%%%THEOREM A PRIORI ESTIMATE%%%%%%%%%% 
\bth {a-prior} There exists two constants $\ga>0$ and $c_{7}>0$ depending on $N$, $p$, $C_{0}$ and $N$, $p$, $C_{0}$ and $R_{0}$ respectively such that if $u$ is any nonnegative solution of (\ref{main}) vanishing on $\prt\Gw\setminus\{0\}$ there holds
\begin {equation}\label {a-prior0}
\myfrac {1}{c_{7}}\gr(x)\abs x^{\ga-1}u(A)\leq
u(x)\leq c_{7} \gr(x)\abs x^{-\ga-1}u(A)
\end {equation}
for any $x\in\Gw$, where $A$ is a fixed point in $\Gw$ such that $\gr(A)\geq R_{0}$.
\es
%%%%%%%%%%%%%%%%%%%%%%%%%%%%%%%%%%%%%%%%%%%%%%%%%%%%%%%%%%%%%%%%%%%%%%
%%%PROOF OF THEOREM A PRIORI ESTIMATE%%%%%
\Proof {\it Step 1}: Tangential estimate. Let $x\in \Gw$ such that
$\abs x=2r\leq R_{0}$ and $\gr(x)=t<br/2 $.
Let $Q\in\prt\Gw\setminus \{0\}$ such that $\abs Q=\abs x$ and
$x\in B_{r}(Q)$, the previous lemma implies
\begin {equation}\label {a-prior1}
\myfrac {2}{c_{6}\abs x}\gr(x)u(N_{r/2}(Q))\leq
u(x)\leq \myfrac {2c_{6}}{\abs x} \gr(x)u(N_{r/2}(Q)).
\end {equation}
There exists a fixed integer $k>2$ such that we can connect two points lying on $\prt B_{2r}(0)\cap\prt\Gw$ by $k$ connected balls $B_{j}$ ($j=1,...,k$) with radius $r/4$ and center on $\prt B_{2r}(0)$. In particular we can connect $N_{r/2}(Q)$ with $N_{2r}(0)=-2r\gn_{0}$ and all the balls can be taken such that the distance of their center to $\prt\Gw$ be larger that $r/2$. Since by \rlemma 
{holdest} there holds
$$
\sup_{B_{j}}u\leq
c_{1}\inf_{B_{j}} u\forevery j=1,...k,
$$
we derive 
\begin {equation}\label {a-prior2}
\myfrac {2}{c_{1}^kc_{6}\abs x}\gr(x)u(N_{2r}(0))\leq
u(x)\leq \myfrac {2c_{1}^kc_{6}}{\abs x} \gr(x)u(N_{2r}(0)).
\end {equation}
Let $A_{0}=-R_{0}\gn_{0}$, $b_{1}=-2r\gn_{0}=N_{2r}(0)$, 
for $\ell\geq 2$, $b_{\ell}=-2(1+3(2^{\ell-1}-1)/2)r\gn_{0}$ and
$r_{\ell}=2^{\ell-1}r$. Applying again \rlemma {holdest} in $B_{2r_{\ell}}(b_{\ell})\subset\Gw$, we have
\begin {equation}\label {a-prior3}
\sup_{B_{r_{\ell}}(b_{\ell})}u\leq
c_{1}\inf_{B_{r_{\ell}}(b_{\ell})} u\forevery \ell=1, 2,...
\end {equation}
Let $\gt$ be the solution of
$$2(1+3(2^{\gt-1}-1)/2)r=R_{0}\Longleftrightarrow
\gt=\myfrac {\ln (R_{0}+r)-\ln 3r}{\ln 2}+1
$$
and $\ell_{0}=\BBE(\gt)+1$, then $A_{0}\in B_{r_{\ell_{0}}}(b_{\ell_{0}})$, and the combination of (\ref {a-prior2}) and
(\ref {a-prior3}) (applied $\ell_{0}$ times) yields to
\begin {equation}\label {a-prior4}
\myfrac {1}{c_{1}^{k+\ell_{0}}c_{6}\abs x}\gr(x)u(A_{0})\leq
u(x)\leq \myfrac {c_{1}^{k+\ell_{0}}c_{6}}{\abs x} \gr(x)u(A_0).
\end {equation}
Since $r\leq R_{0}/2$, the computation of $\gt$ yields to
$$2^{\gt}=\myfrac {2(R_{0}+r)}{3r}\leq \myfrac {R_{0}}{r}\Longrightarrow
c_{1}^\gt\leq \left(\myfrac {R_{0}}{r}\right)^{\ln c_{1}/\ln 2}.
$$
This implies (\ref {a-prior0}) with $\ga=\ln c_{1}/\ln 2$. \medskip

\noindent {\it Step 2}: Internal estimate. If $x\in\Gw$ satisfies
$\abs x\leq R_{0}$ and $\gr (x)\geq b/4\ abs x$, we can directly procede without using \rlemma {two-side}. Using internal Harnack inequality 
(\ref {harn-int}) and connecting $x$ to $N_{r}(0)$ and then to $A_{0}$ we obtain
\begin {equation}\label {a-prior5}
\myfrac {1}{c_{7}}\abs x^{\ga}u(A)\leq
u(x)\leq c_{7} \abs x^{-\ga}u(A),
\end {equation}
from which (\ref {a-prior1}) is derived since $\gr(x) \geq b\abs x/4$. Finally, if $\abs x\geq R_{0}$ and $\gr(x)\leq R_{0}$, we can replace the singular point $0$ by a regular point $B\in\prt\Gw$ such that 
$\abs {x-B}\leq R_{0}$. The previous procedure leads to the same estimate. At end, if $\gr(x)>R_{0}$ we apply again the internal Harnack inequality (\ref{harn-int}). Since $\Gw$ is bounded, $x$ and $A_{0}$ can be joined by at most $d=2\,$diam$\,(\Gw)/R_{0}$ balls  $B_{i}$ with radius $R_{0}/2$ and center $b_{i}$ satisfying
$\gr (b_{i})\geq R_{0}$. Then using $d$ times (\ref {harn-int}) yields to 
(\ref {a-prior1}).\qeda \\

%%%%%%%%%%%%%%%%%%%%%%%%%%%%%%%%%%%%%%%%%%%%%%%%%%%%%%%%%%%%%%%%%%%%%%
%%%REMARK ON THE CONICAL GEOMETRY%%%%%
\noindent \Remark If $p=2$ and $d$ is regular, it is proved that \rlemma {two-side} holds even if $\prt\Gw$ is Lipschitz continuous. The previous proof is adapted by replacing the doubling property of the radius on the connecting balls $B_{r_{\ell}}(b_{\ell})$ by radii such that
$r_{\ell+1}=\gb r_{\ell}$, where $\gb>0$ depends on the opening of the standard cone $C$ associated to the inside cone property of $\Gw$. This observation shows that in the general case $p\neq 2$ and $d$ singular, the validity of \rlemma  {two-side} implies \rth {a-prior} when 
$\Gw$ is a bounded domain satisfying the inside cone property.\\

%%%%%%%%%%%%%%%%%%%%%%%%%%%%%%%%%%%%%%%%%%%%%%%%%%%%%%%%%%%%%%%%%%%%%%
%%%REMARK ON THE Bound%%%%%
\noindent \Remark In the case $p=2$ and $\lim_{x\to 0}\abs {x}^2d(x)=0$ the value of $\ga$ is known and equal to $N-1$. When $p\neq 2$ the value of $\ga$ is unknown, even in the case where 
$d=0$.\\

The next result is a consequence of the method used in the proof of \rth {a-prior}.
%%%%%%%%%%%%%%%%%%%%%%%%%%%%%%%%%%%%%%%%%%THEOREM UNIFORM%%%%%%%%%%
%%%%%%%%%%%%%%%%%%%%%%%%%%%%%%
\bth {unif} Let $u\in C^1(\bar\Gw\setminus \{0\})$ be a positive solutions of (\ref{main}) vanishing on $B_{2R_{0}}\cap(\prt\Gw\setminus \{0\})$. Then there exists a constant $c_{9}>0$ depending on $p$, $N$, $C_{0}$ and $R_{0}$ such that 
\begin {equation}\label {unif1}
\myfrac {1}{c_{9}}\myfrac {u(y)}{\gr(y)}\leq \myfrac {u(x)}{\gr(x)}
\leq c_{9}\myfrac {u(y)}{\gr(y)},
\end {equation}
for every $x$ and $y$ in $B_{R_{0}}(0)\cap\Gw$ satisfying
$\abs y/2\leq\abs x\leq 2\abs y$. 
\es
%%%%%%%%%%%%%%%%%%%%%%%%%%%%%%%%%%%%%%PROOF OF %THEOREM UNIFORM%%%%%%%%%%
%%%%%%%%%%%%%%%%%%%%%%%%%%%%%%
\Proof By (\ref {a-prior2}) we have 
$$
\myfrac {1}{c'\abs x}u(N_{\abs x}(0))\leq
\myfrac {u(x)}{\gr(x)}\leq \myfrac {c'}{\abs x}u(N_{\abs x}(0)).
$$
for any $x\in\Gw$ such that $\abs x\leq R_{0}/2$ and $\gr(x)\leq b\abs x/4$. If we assume that $x\in\Gw\cap B_{R_{0}/2}(0)$ verifies $\abs x\leq R_{0}/2$ and $\gr(x)> b\abs x/4$ we can connect $x$ to 
$N_{\abs x}(0)$ by a fixed number $n$ of balls of radius $b\abs x/8$ with their center at a distance to $\prt\Gw$ larger than $b\abs x$. The classical Harnack inequality yields to
$$
\myfrac {1}{c^n_{1}}u(N_{\abs x}(0))\leq
u(x)\leq c^n_{1}u(N_{\abs x}(0)).
$$
Since $\gr(x)\leq \abs x\leq \gr(x)/b$, we obtain, for any $x\in B_{R_{0}/2}(0)\cap\Gw$,
\begin {equation}\label {unif2}
\myfrac {1}{c_{8}\abs x}\gr(x)u(N_{\abs x}(0))\leq
u(x)\leq \myfrac {c_{8}}{\abs x} \gr(x)u(N_{\abs x}(0)).
\end {equation}
where $c_{8}$ depends on $p$, $N$, $C_{0}$ and $R_{0}$. By Harnack inequality, we can replace $u(N_{\abs x}(0))$ by $u(N_{s}(0))$ for any $\abs x/2\leq s\leq 2\abs x$ and get
\begin {equation}\label {unif3}
\myfrac {1}{c_{1}c_{8}\abs x}\gr(x)u(N_{s}(0))\leq
u(x)\leq \myfrac {c_{1}c_{8}}{\abs x} \gr(x)u(N_{s}(0)).
\end {equation}
If 
$y\in B_{R_{0}/2}(0)\cap\Gw$ satisfies $\abs x/2\leq\abs y\leq \abs x$, 
we apply twice (\ref {unif3}) and we get (\ref {unif2}) with $c_{9}=c^2_{1}c_{8}^2$.\qeda\\

Another consequence of this method and of \rlemma {holdest} and \rlemma {norm} is the 
%%%%%%%%%%%%%%%%%%%%%%%%%%%%%%%
%%%%%THEOREM VERTICAL%%%%%%%%%%%%%%%
%%%%%%%%%%%%%%%%%%%%%%%%%%%%%%%
\bth {vertic}There exists a constant $c'_{9}$ depending on $N$, $p$, $C_{0}$ and $R_{0}$ such that any $u\in C^1(\bar\Gw\setminus \{0\})$ be a positive solutions of (\ref{main}) vanishing on $B_{2R_{0}}\cap(\prt\Gw\setminus \{0\})$ verifies 
\begin {equation}\label {unif1'}
u(x)\leq c'_{9}u(N_{r}(0))
\end {equation}
for every $0<r\leq R_{0}/2$ and any $x\in \Gw\cap B_{2r}(0)\setminus B_{r/2}(0)$.
\es

%%%%%%%%%%%%%%%%%%%%%%%%%%%%%%%
%%%%%REMARK%%%%%%%%%%%%%%%
%%%%%%%%%%%%%%%%%%%%%%%%%%%%%%%
\noindent\Remark Since \rlemma {holdest} and \rlemma {norm} are valid in Lipschitz continuous domains and the construction of connected chain of balls too by \rlemma {geom1}, the above inequality remains valid if $\Gw$ is Lipschitz continuous.\\
%%%%%%%%%%%%%%%%%%%%%%%%%%%%%%%%%%%%%%%%%%%%%%%%%%%%%%%%%%%%%%%%%%%%%%
%%%THEOREM BOUNDARY HARNACK INEQUALITY 

The next result is known as the boundary Harnack inequality.

\bth {bhi} Let $Q\in\prt\Gw$, $0<r\leq \min\{R_{0}/2,\abs Q/2\}$, and $u_{1}$ and $u_{2}$ 
be two nonnegative solutions of (\ref {main}) in $B_{2r}(Q)\cap\Gw$ which vanish on $B_{2r}(Q)\cap\prt\Gw$. Then there exists $c_{10}>0$ depending respectively on $N$, $p$ and $C_{0}$ such that 
\begin {equation}\label {bhi1}
\myfrac {1}{c_{10}}\myfrac {u_{1}(x)}{u_{1}(y)}\leq
\myfrac {u_{2}(x)}{u_{2}(y)}\leq c_{10} \myfrac {u_{1}(x)}{u_{1}(y)}
\end {equation}
for any $x,y\in B_{r}(Q)\cap\Gw$.
\es
%%%%%%%%%%%%%%%%%%%%%%%%%%%%%%%%%%%
%%%PROOF OF THEOREM BOUNDARY HARNACK INEQUALITY %%%%%%%%%%%%%%%%%%%%%%%%%%%%%%%%

\Proof If $x\in B_{r}(Q)\cap\Gw$ satisfies $\gr(x)\leq br/2$, we denote by $P_{x}=P$ the unique projection of $x$ on $\prt\Gw$ and put $t=\gr(x)$. By (\ref{norm-est2}),
\begin {equation}\label {norm-est9}
\myfrac {t}{c_{6}r}\leq \myfrac {u_{i}(x)}{u_{i}(N_{r/2}(Q))}\leq\myfrac {c_{6}t}{r}
\end {equation}
\begin {equation}\label {norm-est10}
\myfrac {1}{c^2_{6}}\myfrac {u_{1}(x)}{u_{1}(N_{r/2}(Q))}\leq
\myfrac {u_{2}(x)}{u_{2}(N_{r/2}(Q))}\leq c^2_{6} \myfrac {u_{1}(x)}{u_{1}(N_{r/2}(Q))}
\end {equation}
from which (\ref {bhi1}) is derived with a first constant $c_{10}=c^4_{6}$. Next, if
$x\in B_{r}(Q)\cap\Gw$ satisfies $\gr(x)> br/2$, we denote $\gb=2+\BBE (-\ln b/\ln 2)$, thus $\gr(x)> r/2^\gb$. By \rlemma {holdestII}
\begin {equation}\label {norm-est11}
 \myfrac {1}{c_2^\gb}\leq \myfrac {u_{i}(x)}{u_{i}(N_{r/2}(Q))}\leq c_2^\gb,
\end {equation}
for $i=1,2$. Therefore (\ref {norm-est10}) holds with $c^2_{6}$ replaced by $c_2^\gb$. Finally (\ref {bhi1}) is verified with $c_{10}=\max\{c^4_{6},c_2^{2\gb}\}$.\qeda \\

The next result is another form of the boundary Harnack inequality
%%%%%%%%%%%%%%%%%%%%%%%%%%%%%%%%%%%
%%%VARIANT OF BOUNDARY HARNACK INEQUALITY  II %%%%%%%%%%%%%%%%%%%%%%%
\bth {bhiII} Let $u_{i}\in C^1(\bar\Gw\setminus \{0\})$ ($i=1,2$) be two nonnegative solutions of (\ref {main}) vanishing on $B_{2R_{0}}\cap(\prt\Gw\setminus \{0\})$. Then there exists $c_{11}>0$ depending respectively on $N$, $p$ and $C_{0}$ such that for any $r\leq R_{0}$
\begin {equation}\label {bhi2}
\sup\left(\myfrac {u_{1}(x)}{u_{2}(x)}:x\in \Gw\cap (B_{r}(0)\setminus B_{r/2}(0)\right)\leq
c_{11}\inf\left(\myfrac {u_{1}(x)}{u_{2}(x)}:x\in \Gw\cap (B_{r}(0)\setminus B_{r/2}(0)\right).
\end {equation}
\es
%%%%%%%%%%%%%%%%%%%%%%%%%%%%%%%%%%%
%%%PROOF%%VARIANT OFBOUNDARY HARNACK INEQUALITY  II %%%%%%%%%%%%%%%%%%%%%%%
\Proof Applying twice \rth {unif}, we get
\begin {equation}\label {bhi3}
\myfrac {1}{c^2_{9}}\myfrac {u_{1}(x)}{u_{1}(y)}\leq \myfrac {u_{2}(x)}{u_{2}(y)}
\leq c^2_{9}\myfrac {u_{1}(x)}{u_{1}(y)},
\end {equation}
for any $x$ and $y$ such that $\abs x/2\leq \abs y\leq 2\abs x$. Equivalently
\begin {equation}\label {bhi4}
\myfrac {1}{c^2_{9}}\myfrac {u_{1}(x)}{u_{2}(x)}\leq \myfrac {u_{1}(y)}{u_{2}(y)}
\leq c^2_{9}\myfrac {u_{1}(x)}{u_{2}(x)},
\end {equation}
which the claim with $c_{11}=c_{9}^2$.\qeda 

\subsection {The set of singular solutions}

We still assume that $\Gw$ is a bounded domain with a $C^2$ boundary containing the singular point $0$. We introduce the following assumption on the function $d$.
%%%%%%%%%%%%%%%%%%%%%%%%%%%%%%%%%%%
%%%DEFINITION LOCAL COMPPRINCIPLE%%%%% %%%%%%%%%%%%%%%%%%%%%%
\bdef {mp} A measurable function $d$ satisfying (\ref {main-d}) with 
$a=0\in\prt\Gw$ is said to satisfy the local comparison principle in $\Gw$ if, for any $\ge>0$ and any $u_{i}\in C^1(\bar\Gw_{\ge})$ ($i=1,2$) nonnegative solutions of (\ref{main}) in $\Gw_{\ge}=\Gw\setminus \bar B_{e}(0)$ which vanish on $\prt^*\Gw_{\ge}=\prt\Gw\setminus B_{\ge}(0)$,  $u_{1}(x)\geq u_{2}(x)$ on
$\Gw\cap \prt B_{\ge}(0)$ implies
 $u_{1}\geq u_{2}$ in $\bar\Gw_{\ge}$.
\es

 Clearly, if $d$ is nonpositive it satisfies the  local comparison principle. However there are many other cases, depending either on the value of $C_{0}$ or the rate of blow-up of $d$ near $0$ which insure this principle. 
 %%%%%%%%%%%%%%%%%%%%%%%%%%%%%%%%%%%
%%%THEOREM UPPER ESTIMATE%%%%%%%%% %%%%%%%%%%%%%%%%%%%%%%
  \bth {upp-est} Assume $d$ satisfies the local comparison principle and there exists a nonnegative nonzero solution $u$ to (\ref {main}) in $\Gw$ which vanishes on $\prt\Gw\setminus \{0\}$. If $v$ is any other nonnegative solution of (\ref {main}) in $\Gw$ vanishing on $\prt\Gw\setminus \{0\}$ there exists $k\geq 0$ such that $v\leq ku$.
\es
%%%%%%%%%%%%%%%%%%%%%%%%%%%%%%%%%%%
%%%PROOF OF THEOREM UPPER ESTIMATE%%%%%%%%% %%%%%%%%%%%%%%%%%%%%%%
\Proof Since any nontrivial nonnegative solution is positive by Harnack inequalities we can assume that both $u$ and $v$ are positive in $\Gw$. We denote by $\CH$ the set of $h>0$ such that 
 $v< hu$ in $\Gw$ and we assume that $\CH$ is empty otherwhile the results is proved. Then for any 
 $n\in \BBN_{*}$ there exists $x_{n}\in\Gw$ such that $v(x_{n})\geq n u(x_{n})$. We can assume that 
 $x_{n}\to\xi$ for some $\xi\in \bar\Gw$. Clearly $\xi\in\Gw$ is impossible. Let us assume first that 
 $\xi\in \prt\Gw\setminus \{0\}$ and denote by $\xi_{n}$ the projection of $x_{n}$ onto $\prt\Gw$. Thus
 $$\myfrac {v(x_{n})-v(\xi_{n})}{\gr (x_{n})}\geq n\myfrac {u(x_{n})-u(\xi_{n})}{\gr (x_{n})}.
 $$
 Because $u$ and $v$ are $C^1$ in $\bar\Gw\setminus \{0\}$, 
 $$\lim_{n\to\infty}\myfrac {v(x_{n})-v(\xi_{n})}{\gr (x_{n})}=\myfrac {\prt v}{\prt\gn_{\xi}}(\xi)
\;\mbox { and }\;
\lim_{n\to\infty}\myfrac {u(x_{n})-u(\xi_{n})}{\gr (x_{n})}=\myfrac {\prt u}{\prt\gn_{\xi}}(\xi).
 $$
 Since Hopf boundary lemma  is valid (see \cite {To}), the two above normal derivative at $\xi$ are negative,
 which leads to a contradiction. Thus we are left with the case $x_{n}\to 0$. Set $r_{n}=\abs {x_{n}}$. By \rth
 {bhiII}
 $$\inf\left\{\myfrac {v(x)}{u(x)}:\abs x=r_{n}\right\}\geq c^{-1}_{11}\myfrac {v(x_{n})}{u(x_{n})}\geq c^{-1}_{11}n
 $$
By the  local comparison principle assumption, $v\geq nc^{-1}_{11}u$ in $\Gw_{r_{n}}$. This again leads to a contradiction.  \qeda
\medskip

 The next statement is useful to characterize unbounded solutions
 %%%%%%%%%%%%%%%%%%%%%%%%%%%%%%%%%%%
%%%PROPOSITION SINGULAR SOLUTION%%%%%%%%% %%%%%%%%%%%%%%%%%%%%%%
 \bprop{sing} Assume $u$ is a nonnegative solution of (\ref {main}) vanishing on $\prt\Gw\setminus\{0\}$,  unbounded and without extremal points near $0$. Then 
  \begin {equation}\label {sing1}
 \lim_{x\to 0}\myfrac{\abs xu(x)}{\gr(x)}=\infty.
 \end {equation}
 \es
 %%%%%%%%%%%%%%%%%%%%%%%%%%%%%%%%%%%
%%%PROOF OF PROPOSITION SINGULAR SOLUTION%%%% %%%%%%%%%%%%% %%%%%%%%%%%%%%%%%%%%%%
 \Proof  Assume that (\ref {sing1}) is not true. Then there exist a sequence $\{s_{n}\}$ converging to $0$ and a constant $M>0$ such that
$$\sup\left\{ \myfrac{\abs xu(x)}{\gr(x)}:\abs x=s_{n}\right\}\leq M.
$$
 Therefore $\sup\left\{u(x):\abs x=s_{n}\right\}\leq M $. Because $u$ has no extremal points near $0$, say in $B_{s_{0}}(0)$ for some $s_{0}>0$, the maximum of $u$ in 
 $\Gw\cap (B_{s_{0}}(0)\setminus B_{s_{n}}(0))$ is achieved either on 
 $\abs x=s_{0}$ or on  $\abs x=s_{n}$. Therefore
 $$\max\{u(x):x\in \Gw\cap (B_{s_{0}}(0)\setminus B_{s_{n}}(0))\}
 \leq \max\left\{M, \max\{u(x):x\in \Gw\cap \prt B_{s_{0}}(0)\}\right\}=M.
 $$
Since this is valid for any $n$, it implies that $u$ is bounded in $\Gw$, contradiction.\qeda \\

Such a solution is called a {\it singular solution}. The next result, which extends a previous result in \cite {Bo}, made more precise the statement of \rth {upp-est}.
 %%%%%%%%%%%%%%%%%%%%%%%%%%%%%%%%%%%
%%%THEOREM 1 DIMENSIONAL SET%%%%%%%%% %%%%%%%%%%%%%%%%%%%%%%
 \bth {1-D} Assume $d$ satisfies the local comparison principle and there exists a positive singular solution $u$ to (\ref {main}) in $\Gw$ vanishing on $\prt\Gw\setminus \{0\}$. Assume also either $1<p\leq 2$ and $d\geq 0$, or 
 $p>2$,  $u$ admits no critical point in $\Gw$ and
  \begin {equation}\label {sing1'}
 \liminf_{x\to 0}\myfrac{\abs x\abs {Du(x)}}{u(x)}>0.
 \end {equation}
If $v$ is any other positive solution of (\ref {main}) in $\Gw$ vanishing on $\prt\Gw\setminus \{0\}$ there exist $k\geq 0$ such that $v=ku$.
\es
 %%%%%%%%%%%%%%%%%%%%%%%%%%%%%%%%%%%
%%%PROOF OF THEOREM 1 DIMENSIONAL SET%%%%%%%%% %%%%%%%%%%%%%%%%%%%%%%
 \Proof Let us assume that $v$ is not zero. By \rth {upp-est} there exists a minimal $k>0$ such that $v\leq ku$. As in the proof of \rth {upp-est} the following holds:\smallskip
 
 \noindent (i) either the graphs of $v$ and $ku$ are tangent at some $\xi\in \Gw$. If we set $w=ku-v$, then $w(\xi)=0$, and
\begin {equation}\label {linea}
-\CL w-Dw=0
\end {equation}
where $\CL$ is a linear elliptic operator and $D=d(x)(k^{p-1}u^{p-1}-v^{p-1})/w$. Since $ku(\xi)=v(\xi)>0$, $D$ is locally bounded near $\xi$.\;If $p>2$ and $u$ admits no critical point in $\Gw$, $\CL$ is uniformly elliptic (\cite {FV}, \cite {Bo} for details in a similar situation). Thus the strong maximum principles holds and $w$ is locally zero. Since $\Gw$ is connected $w\equiv 0$ in $\Gw$. If $1<p\leq 2$, the strong maximum principle holds to and we have the same conclusion.\smallskip
 
 \noindent (ii) either the graphs of $v$ and $ku$ are not tangent inside $\Gw$, but tangent on $\prt \Gw\setminus \{0\}$. Since the normal derivatives of $ku$ and $v$ at $\xi$ coincide, $\CL$ is uniformly elliptic. If $p\geq 2$ the coefficient $D$ is locally bounded. If $1<p\leq 2$ this is not the case but $D$ remains nonnegative. In both case Hopf maximum principle applies and yields to $\prt w/\prt\gn_{\xi}(\xi)<0$. This is again a contradiction.
 \smallskip
 
 \noindent (iii) or $v< ku$ in $\Gw$, $\prt v/\prt\gn>k\prt u/\prt\gn$ on 
 $\prt \Gw\setminus \{0\}$ and there exists a sequence $\{x_{n}\}\subset \Gw$ converging to $0$ such that
 $$\lim_{n\to\infty}\myfrac {v(x_{n})}{u(x_{n})}=k.
 $$
Furthermore we can assume that 
$$\myfrac {v(x_{n})}{u(x_{n})}=\sup\left\{\myfrac {v(x)}{u(x)}:
\abs x=\abs {x_{n}}=:r_{n}\right\}
$$
Put $a_{n}=\max\{u(x):\abs x=r_{n}\}$. By \rth {vertic} there exists $c'_{9}>0$ depending on $N$, $p$, $C_{0}$ and $R_{0}$ such that
\begin {equation}\label {high1}
u(N_{r_{n}}(0))\leq a_{n}\leq c'_{9}u(N_{r_{n}}(0)),
\end {equation}
which implies
\begin {equation}\label {high2}
\max\{u(x):r_{n}/2\leq \abs x\leq 2r_{n}\}\leq c'_{9}a_{n}.
\end {equation}
We set
 $u_{n}(x)=u(r_{n}x)/a_{n}$, $v_{n}(x)=v(r_{n}x)/a_{n}$ and $d_{n}(x)=r_{n}^pd(r_{n}x)$. Then both $u_{n}$ and $v_{n}$ are solutions of
 $$-div\left(\abs {Df}^{p-2}Df\right)-d_{n}f^{p-1}=0
 $$
 in $\Gw_{n}=\Gw/r_{n}$ and vanish on $\prt\Gw_{n}\setminus\{0\}$. By (\ref {high2}), $u_{n}$ and $v_{n}$ are uniformly bounded in 
 $\tilde\Gw_{n}=\Gw_{n}\cap (B_{2}(0)\setminus\bar B_{1/2}(0))$. Since 
$\prt \Gw_{n}\cap(B_{2}(0)\setminus\bar B_{1/2}(0)) $ is uniformly 
$C^2$ we deduce by the degenerate elliptic equations theory \cite {Lie} that, up to subsequences, $u_{n}$ and $v_{n}$ converge in the 
 $C^1_{loc}(\bar\Gw_{n}\cap (B_{2}(0)\setminus\bar B_{1/2}(0)))$-topology to functions $U$ and $V$ which satisfy
  $$-div\left(\abs {Df}^{p-2}Df\right)-d_{\infty}f^{p-1}=0
 $$
 in $H\cap (B_{2}(0)\setminus\bar B_{1/2}(0)))$, where $H$ is the half space $\{\eta\in \BBR^N:\eta.\gn_{0}<0\}$ and $d_{\infty}$ is some weak limit of $d_{n}$ in the
 weak-star topology of $L^\infty$. Moreover, if $p>2$, (\ref {sing1'}), jointly with (\ref{high1}) and (\ref{high2}), implies that
\begin {equation}\label {sing1''}
\abs {Du_{n}(x)}=\myfrac {r_{n}Du(r_{n}x)}{a_{n}}\geq \gg
\end {equation}
for $2/5\leq \abs x\leq 8/5$, where $\gg>0$.
 We put
 $$
 \ell_{n}=\inf\left\{\myfrac {v(x)}{u(x)}:
\abs x=\abs {x_{n}}=:r_{n}\right\}\leq k
$$
 and $\xi_{n}=x_{n}/r_{n}$. Up to another choice of subsequence, we can also assume that $\ell_{n}\to\ell$ and $\xi_{n}\to\xi$ with $\xi=1$. Furthermore $V\leq kU$, $V(\xi)=kU(\xi)$ and, if 
 $\xi\in \prt H\cap (B_{2}(0)\setminus\bar B_{1/2}(0)))$,
 $$\myfrac {\prt V}{\prt \gn_{0}}(\xi)=k\myfrac {\prt U}{\prt \gn_{0}}(\xi) <0.
 $$
In this case, and more generally if the coincidence set $\Xi$ of $V$ and $kU$ has a nonempty intersection with  $\prt H\cap (B_{2}(0)\setminus\bar B_{1/2}(0)))$, Hopf boundary lemma applies and implies that $V=kV$ in the whole domain. If this is not the case we use (\ref {sing1''}) to
conclude again by the strong maximum principle that $V=kU$ in $H\cap (B_{2}(0)\setminus\bar B_{1/2}(0)))$. Therefore $\ell=k$ and for any $\ge>0$ there exists $n_{\ge}\in \BBN$ such that $n\geq n_{\ge}$ implies
$$(k-\ge)u(x)\leq v(x)\leq ku(x)\forevery x\in\Gw\cap \prt B_{r_{n}}(0).
$$
By the local comparison principle the same estimate holds in $\Gw_{r_{n}}$. Since this is valid for any $n$ and any $\ge$, we conclude that $v=ku$.
\qeda
%%%%%%%%%%%%%%%%%%%%%%%%%%%%%%%%%%%%%%%%%%%%%%%%%%%%%%
%%%%%%%%EXISTENCE OF SINGULAR SOLUTIONS%%%%%%%%%%%%%%%%%%%%%%%
%%%%%%%%%%%%%%%%%%%%%%%%%%%%%%%%%%%%%%%%%%%%%%%%%%%%%%
 \mysection {Existence of singular solutions}
 %%%%%%%%%%%%%%%%%%%%%%%%%%%%%%%%%%%%%%%%%%%%%%%%%%%%%%
%%%%%%%%SEPARABLE SOLUTIONS%%%%%%%%%%%%%%%%%%%%%%%
%%%%%%%%%%%%%%%%%%%%%%%%%%%%%%%%%%%%%%%%%%%%%%%%%%%%%%

\subsection {Separable solutions}
The existence of N-dimensional regular separable p-harmonic functions associated to cones is due to Tolksdorff \cite {To}. Extension to singular function is proved in 
\cite {Ve2}. These solutions are obtained as follows: {\it Let $(r,\gs)\in \BBR_{+}\ti S^{N-1}$ be the spherical coordinates in $\BBR^N$ and $S\subset S^{N-1}$ a smooth spherical domain. Then there exists two couples $(\gg_{S}, \psi_{S})$ and 
 $(\gb_{S}, \phi_{S})$, where $\gg_{S}$ and $\gb_{S}$ are positive real numbers and 
$\psi_{S}$ and $\phi_{S}$ belong to $C^2(\bar S)$ and vanish on $\prt S$,  such that 
\begin {equation}\label {spheric}
U_{S}=r^{\gg_{S}}\psi_{S}\quad\mbox {and }\;V_{S}=r^{-\gb_{S}}\phi_{S},
\end {equation}
are p-harmonic functions in the cone $C_{S}=\{(r,\gs):r>0,\,\gs\in S\}$. These couples are unique up to homothety over $\psi_{S}$ and $\phi_{S}$. Furthermore
the following equation holds 
\begin {equation}\label {spheric1}\left\{ \BA {l}
-div\left((a^{2}\eta^2+\abs {\nabla\eta}^2)^{(p-2)/2}\nabla\eta\right)
=\gl(a)(a^{2}\eta^2+\abs {\nabla\eta}^2)^{(p-2)/2}\eta\quad\mbox {in }S\\
[2mm]\phantom {-div,,(a^{2}\eta^2+\abs {\nabla\eta}^2)^{(p-2)/2}\nabla\eta}
\eta=0\quad\mbox {on } \prt S.
\EA\right.\end {equation}
where 
$\gl(a)=a(a(p-1)+p-N)$ if $(a,\eta)=(\gb_{S},\phi_{S})$, and 
$\gl(a)=a(a(p-1)+N-p)$ if $(a,\eta)=(\gg_{S},\psi_{S})$.}\\

If $p\neq 2$ and $N\neq 2$,  $\gg_{S}$ and $\gb_{S}$ are unknown except if 
$S=S^{N-1}_{+}=S^{N-1}\cap \{x=(x',x_{N}):x_{N}>0\in\BBR^N\}$, in which case $\gg_{S}=1$ and 
$\psi_{S}=x_{N}$. If $N=2$, equation (\ref{spheric1}) is completely integrable and the values of the $\gg_{S}$ and $\gb_{S}$ are known (\cite {KM}, \cite {KV}). When $p\neq 2$, the existence of solutions to (\ref {spheric1}) is not easy since this is not a variational problem on $S$. 
 Tolksdorff's method (\cite {To}) is based upon a $N$-dimensional shooting argument: he constructs the solution
 $v$ of
 \begin{equation} \label {p-harm}\left\{\BA {l}
 -div\left(\abs {Dv}^{p-2}Dv\right)=0\quad\mbox {in }\;C^1_{S}=C_{S}\cap\{x:\abs x\geq 1\}\\
\phantom {--------,} v=(2-\abs x)_{+}\quad\mbox {on }\;\prt C^1_{S}.
 \EA\right.\end {equation}
 Then he proves, thanks to an equivalence principle,  that the function $v$ stabilizes at infinity under the asymptotic form $v(x)\approx \abs x^{-\gb}\gf (x/\abs x)$, with $\gb>0$, which gives (\ref {spheric1}) and the function $V_{S}$. The domain $S$ characterizes the exponent $\gb$. The same argument applies if (\ref {p-harm}) is replaced by
  \begin{equation} \label {p-harm2}\left\{\BA {l}
 -div\left(\abs {Dv}^{p-2}Dv\right)+\myfrac {c u^{p-1}}{\abs x^p}=0\quad\mbox {in }\;C^1_{S}=C_{S}\cap\{x:\abs x\geq 1\}\\
\phantom {----,--\myfrac {c u^{p-1}}{\abs x^p}----,} 
v=(2-\abs x)_{+}\quad\mbox {on }\;\prt C^1_{S}.
 \EA\right.\end {equation}
 with $c>0$. This gives rise to a solution of (\ref {p-harm2}) in $C_{S}$ under the form $V_{S,c}=r^{-\gb_{c,S}}\eta$ where
 $\gb_{c,S}>0$ and 
 \begin {equation}\label {spheric3}\left\{ \BA {l}
-div\left((\gb_{c,S}^{2}\eta^2+\abs {\nabla\eta}^2)^{(p-2)/2}\nabla\eta\right)+c\eta^{p-1}
=\gl(\gb_{c,S})(\gb_{c,S}^{2}\eta^2+\abs {\nabla\eta}^2)^{(p-2)/2}\eta\quad\mbox {in }S\\
[2mm]\phantom {-div,,(a^{2}\eta^2+\abs {\nabla\eta}^2)^{(p-2)/2}\nabla\eta+c\eta^{p-1},}
\eta=0\quad\mbox {on } \prt S.
\EA\right.\end {equation}
With these considerations we can construct singular solutions of (\ref{main}) under a very restrictive geometry
for $\Gw$, where $S=S_{+}^{N-1}$, the upper half unit sphere.

\bth{singsol} Assume $d(x)=-c\abs x^{-p}$ with $c\geq 0$ and $\Gw$ is a bounded domain with a $C^2$ boundary containing $0$. Assume also $\prt\Gw$ is flat in a neighborhood of $0$ and $x.\gn_{0}\leq 0$ for any 
$x\in\bar\Gw$. Then there exists a positive solution of (\ref{main}) which vanishes on $\prt\Gw\setminus \{0\}$
and satisfies
 \begin {equation}\label {blow-up}
 \lim_{\tiny{\BA {l}x\to 0\\ x/\abs x\to\gs\EA}}
 \abs x^{\gb_{c,S_{+}^{N-1}}}u(x)=\eta(\gs)
 \end {equation}
 uniformly for $\gs\in S_{+}^{N-1}$, where $\eta$ is a positive solution of (\ref{spheric3}).
 \es
\Proof Since $\Gw$ is located on one side of the hyperplane $\Gl=\{x:x.\gn_{0}= 0\}$, the restriction to $\Gw$
of the function $V_{c,S_{+}^{N-1}}:x\mapsto  \abs x^{-\gb_{c,S_{+}^{N-1}}}\eta(x/\abs x)$ is a singular solution of 
 \begin {equation}\label {main-D}
-div\left({\abs {Dv}}^{p-2}Dv\right)+\myfrac {c}{\abs x^p}\abs v^{p-2}v=0
 \end {equation}
which vanishes on $\Gl\cap \prt\Gw\setminus\{0\}$ and is positive on $\bar\Gw\setminus H$.
Let $K=\max\{V_{c,S_{+}^{N-1}}(x):x\in \prt\Gw\setminus \Gl\}$. Then $V_{c,S_{+}^{N-1}}-K$ is a solution (or a subsolution if $c>0$). For any $\ge>0$ let $u_{\ge}$ be the solution of
 \begin {equation}\label  {main-D1}\left\{ \BA {l}
-div\left(\abs {Du_{\ge}}^{p-2}Du_{\ge}\right)+\myfrac {c}{\abs x^p}\abs{u_{\ge}}^{p-2}u_{\ge}=0\quad\mbox {in }\Gw\setminus B_{\ge}(0)\\
[2mm]u_{\ge}=0\quad\mbox {on } \prt\Gw\cap B^c_{\ge}(0)
\\
[2mm]u_{\ge}=V_{c,S_{+}^{N-1}}\quad\mbox {on } \Gw\cap \prt B_{\ge}(0).
\EA\right.\end {equation}
Then $V_{c,S_{+}^{N-1}}-K\leq u_{\ge}\leq V_{c,S_{+}^{N-1}}$ and the mapping $\ge\mapsto u_{\ge}$ is increasing. Therefore $u_{\ge}$ converges to $u$ in the $C^1_{loc}(\bar\Gw\setminus\{0\})$ topology and $u$ is a solution of (\ref {main}) which vanishes on $\prt\Gw\setminus\{0\}$ and satisfies (\ref {blow-up}).
\qeda\\

\noindent \Remark If $N=p$ the set of p-harmonic functions is invariant under the M\"oebius group, and in particular under the transformation $x\mapsto \CI(x)=x/\abs x^2$ which preserves $C_{S}$. In such a case $\gb_{S^{N-1}_{+}}=1$. By using the transformation $\CI$ it is possible to prove (see \cite {Bo}) that there exist positive 
$N$-harmonic functions in any bounded domain $\Gw$ having a singularity at a point $a$ of the boundary 
and vanishing on $\prt\Gw\setminus \{a\}$. \\

\noindent \Remark When $p=2$ it is possible to prove the existence of a singular solution to 
 \begin {equation}\label  {lin1}\left\{ \BA {l}
-\Gd u+d(x)u=0\quad\mbox {in }\Gw\\
u=0\quad\mbox {on } \prt\Gw\setminus\{a\}
\EA\right.\end {equation}
where $a\in \prt\Gw$, for any $C^2$ domain $\Gw$ and any $d$ locally bounded in $\bar\Gw\setminus \{a\}$ such that 
$$-\infty<\liminf_{x\to a}\abs {x-a}^2d(x)\leq\limsup_{x\to a}\abs {x-a}^2d(x)< N^2/4.
$$
We conjecture that such a result holds for (\ref{main}) and $p\neq 2$ although the precise upper limit as $x\to a$ of $\abs {x-a}^pd(x)$. We believe that  at least if $N\geq p$ and
$\limsup_{x\to a}\abs {x-a}^2d(x)\leq ((N-p)/p)^p$, (the Hardy constant for $W^{1,p}$ in $\BBR^N$), such a singular solution do exist.

%%%%%%%%%%%%%%%%%%%%%%%%%%%%%%%%%%%%%%%%%%%
\small{
Laboratoire de Math\'ematiques et Physique Th\'eorique\\
CNRS UMR 6083\\
Facult\'e des Sciences\\
Universit\'e Fran\c{c}ois Rabelais\\
F37200 Tours  France\\

\noindent veronmf@univ-tours.fr\\
borghol@univ-tours.fr\\
veronl@lmpt.univ-tours.fr}
%Furthermore%Furthermore%Furthermore
%Furthermore%Furthermore%Furthermore
%Furthermore%Furthermore%Furthermore
 %%END DOCUMENT%%%%%%%%%%%%%%%%%%%%%%%%
\end {document}